# Faà di Bruno's note on eponymous formula  
# Trilingual version

Julyan Arbel

Inria Grenoble - Rhône-Alpes

December 15, 2016

**Abstract**

About 160 years ago, the Italian mathematician Faà di Bruno published two notes dealing about the now eponymous formula giving the derivative of any order of a composition of two functions. We reproduce here the two original notes, Faà Di Bruno (1855, 1857), written respectively in Italian and in French, and propose a translation in English.

## Introduction

The Italian mathematician Francesco Faà di Bruno was born in Alessandria (Piedmont, Italy) in 1825 and died in Turin in 1888. At the time of his birth, Piedmont used to be part of the Kingdom of Sardinia, led by the Dukes of Savoy. Italy was then unified in 1861, and the Kingdom of Sardinia became the Kingdom of Italy, of which Turin was declared the first capital. At that time, Piedmontese used to commonly speak both Italian and French.

Faà di Bruno is probably best known today for the eponymous formula which generalizes the derivative of a composition of two functions, $\phi \circ \psi$, to any order:

$$(\phi \circ \psi)^{(n)} = \sum \frac{n!}{m_1! \ldots m_n!} \phi^{(m_1 + \cdots + m_n)} \circ \psi \cdot \prod_{i=1}^{n} \left(\frac{\psi^{(j)}}{j!}\right)^{m_j}$$

over $n$-tuples $(m_1, \ldots, m_n)$ satisfying $\sum_{j=1}^{n} jm_j = n$.

Faà di Bruno published his formula in two notes[1] referenced as Faà Di Bruno (1855) in Italian, and Faà Di Bruno (1857) in French. They both date from December 1855, and were signed in Paris. They are similar and essentially state the formula without a proof. For a detailed historical review of the formula, and references which predate Faà di Bruno's works such as Abadie (1850), see Johnson (2002).

The rest of this note[2] contains a translation from the French version to English on page 2, the original version in French on page 4, and the original version in Italian on page 6.

---

[1] Original notes can be accessed online on Google Books.
[2] This note is typeset using the Erasmus MMXVI font.

# NOTE ON A NEW FORMULA FOR DIFFERENTIAL CALCULUS.

By M. Faà di Bruno.

Having observed, when dealing with series development of functions, that there did not exist to date any proper formula dedicated to readily calculate the derivative of any order of a function of function without resorting to computing all preceding derivatives, I thought that it would be very useful to look for it. The formula which I found is well simple; and I hope it shall become of general use henceforth.

Let $\phi(x)$ be any function of the variable $x$, linked to another one $y$ by the equation of the form

$$x = \Psi(y).$$

Denote by $\Pi(l)$ the product[3] $1 \cdot 2 \cdot 3 \ldots l$ and by $\Psi'$, $\Psi''$, $\Psi'''$, etc the successive derivatives of the function $\Psi(y)$; the value of $D_y^n \phi$ will have the following expression:

$$D_y^n \phi = \sum \frac{\Pi(n)}{\Pi(i)\Pi(j)\ldots\Pi(k)} \times$$
$$D_x^p \phi \cdot \left(\frac{\Psi'}{1}\right)^i \left(\frac{\Psi''}{1 \cdot 2}\right)^j \left(\frac{\Psi'''}{1 \cdot 2 \cdot 3}\right)^h \cdots \left(\frac{\Psi^{(l)}}{\Pi(l)}\right)^k.$$

the sign $\sum$ runs over all integer and non negative values of $i, j, h \ldots k$, for which

$$i + 2j + 3h \ldots + lk = n,$$

the value of $p$ being

$$p = i + j + h \ldots + k$$

---

[3] *product* written in English in the French version.

The expression can also take the form of a determinant, and one has

$$(-)^n D_y^{n+1}\phi =$$

$$\begin{vmatrix} \psi^{(n+1)}\phi, & \psi'\phi, & \frac{n}{1}\psi''\phi, & \frac{n\cdot n-1}{1\cdot 2}\psi'''\phi, & \ldots\ldots, & \frac{n}{1}\psi^{(n)}\phi \\ \psi^{(n)}\phi, & -1, & \psi'\phi, & \frac{n-1}{1}\psi''\phi, & \ldots\ldots, & \frac{n-1}{1}\psi^{(n-1)}\phi \\ \psi^{(n-1)}\phi, & 0, & -1, & \psi'\phi, & \ldots\ldots, & \frac{n-2}{1}\psi^{(n-2)}\phi \\ \psi^{(n-2)}\phi, & 0, & 0, & -1, & \ldots\ldots, & \frac{n-3}{1}\psi^{(n-3)}\phi \\ \vdots & & & & & \vdots \\ \psi''\phi, & 0, & 0, & 0, & \ldots -1, & \psi'\phi \\ \psi'\phi, & 0, & 0, & 0, & 0, & -1 \end{vmatrix}$$

It is implicit that the exponents of $\phi$ will be considered as orders of derivation.

*Paris, December* **17. 1855.**

---

# NOTE SUR UNE NOUVELLE FORMULE DE CALCUL DIFFÉRENTIEL.

Par M. Faà di Bruno.

AYANT vu, en m'occupant du développement des fonctions en série, qu'il n'existait jusqu'à présent aucune formule propre à calculer immédiatement la dérivée d'un ordre quelconque d'une fonction de fonction sans passer par le calcul de toutes les dérivées précédentes, j'ai pensé qu'il serait très utile de la chercher. La formule, que j'ai trouvée, est bien simple ; et j'espère qu'elle deviendra désormais d'un emploi général.

Soit $\phi(x)$ une fonction quelconque de la variable $x$, liée à une autre $y$ par une équation de la forme

$$x = \Psi(y).$$

Désignons par $\Pi(l)$ le product[4] $1.2.3\ldots l$ et par $\Psi'$, $\Psi''$, $\Psi'''$, etc. les dérivées successives de la fonction $\Psi(y)$ ; on aura pour la valeur de $D_y^n \phi$ l'expression suivante :

$$D_y^n \phi = \sum \frac{\Pi(n)}{\Pi(i)\Pi(j)\ldots\Pi(k)} \times$$
$$D_x^p \phi \cdot \left(\frac{\Psi'}{1}\right)^i \left(\frac{\Psi''}{1.2}\right)^j \left(\frac{\Psi'''}{1.2.3}\right)^h \cdots \left(\frac{\Psi^{(l)}}{\Pi(l)}\right)^k.$$

le signe $\sum$ s'étendant à toutes les valeurs entières et positives de $i, j, h \ldots k$, pour lesquelles

$$i + 2j + 3h \ldots + lk = n,$$

---

[4] The French word should be *produit* here.

et la valeur de $p$ etant[5]

$$p = i + j + h \ldots + k$$

Cette expression peut se mettre aussi la forme[6] d'un déterminant, et l'on a,

$$(-)^n D_y^{n+1} \phi =$$

$$\begin{vmatrix} \psi^{(n+1)}\phi, & \psi'\phi, & \frac{n}{1}\psi''\phi, & \frac{n \cdot n-1}{1 \cdot 2}\psi'''\phi, & \ldots\ldots, & \frac{n}{1}\psi^{(n)}\phi \\ \psi^{(n)}\phi, & -1, & \psi'\phi, & \frac{n-1}{1}\psi''\phi, & \ldots\ldots, & \frac{n-1}{1}\psi^{(n-1)}\phi \\ \psi^{(n-1)}\phi, & 0, & -1, & \psi'\phi, & \ldots\ldots, & \frac{n-2}{1}\psi^{(n-2)}\phi \\ \psi^{(n-2)}\phi, & 0, & 0, & -1, & \ldots\ldots, & \frac{n-3}{1}\psi^{(n-3)}\phi \\ \vdots & & & & & \vdots \\ \psi''\phi, & 0, & 0, & 0, & \ldots -1, & \psi'\phi \\ \psi'\phi, & 0, & 0, & 0, & 0, & -1 \end{vmatrix}$$

Il est sousentendu que les exposants de $\phi$ devront être considérés comme des ordres de différentiation.

*Paris, Decembre* 17. 1855.

---

[5] étant

[6] one should read 'sous la forme'

# SULLO SVILUPPO DELLE FUNZIONI

*NOTA*

**DEL CAV F. FAA' DI BRUNO.**

Nello sviluppo delle funzioni e nelle applicazioni della serie di Lagrange egli arriva sovente di dover cercare la derivata di un dato ordine di una funzione di funzione. Finora tale ricerca richiedeva il calcolo preliminare delle derivate precedenti, cosa penosissima, allorquando l'ordine è grande, come nella teoria delle perturbazioni. Mi pare pertanto che sarebbe stato cosa ben utile l'avere una formola che somministrasse immediatamente l'espressione di tale derivata. Eccola pertanto come mi venne' fatto di trovare.

Sia $\phi(x)$ una funzione qualunque della variabile $x$ legata ad un altra $y$ per la

$$x = \Psi(y).$$

Si avrà:

$\mathbf{D}_y^n \phi =$

$$\sum \frac{\pi(n)}{\pi i \,.\, \pi j \,.\, \pi h \ldots \pi l} [\mathbf{D}_x^p \phi] y \left(\frac{\psi'}{1}\right)^i \left(\frac{\psi''}{1\,.\,2}\right)^j \left(\frac{\psi'''}{1\,.\,2\,.\,3}\right)^h \ldots \left(\frac{\psi^{(l)}}{1\,.\,2\,.\,3\ldots l}\right)^k$$

ove $\pi(l) = 1\,.\,2\,.\,3\ldots l$ ed il segno $\sum$ estendesi a tutti i valori interi e positivi di $p, i, j, h, \ldots l$ che verificano l'equazione

$$i + 2j + 3h + \ldots + lk = n$$
$$p = i + j + h + \ldots + k$$

Allorquando, come nella serie di Lagrange, si ha:

$$\phi(x) = x^m$$

si avrà sotto le medesime condizioni :

$$\frac{\mathbf{D}^n(\psi y)^m}{1\,.\,2\,.\,3\ldots n} =$$
$$\sum \frac{m(m-1)\ldots(m-p+1)}{\pi i\,.\,\pi j\,.\,\pi h\ldots}\psi^{m-p}\left(\frac{\psi'}{1}\right)^i\left(\frac{\psi''}{1\,.\,2}\right)^j\left(\frac{\psi'''}{1\,.\,2\,.\,3}\right)^h\ldots$$

L'espressione generale di $\mathbf{D}_y^n\phi$ puossi mettere sotto la forma di un determinante. Così a modo d'esempio

$$\begin{vmatrix} \psi'''\cdot\phi & \psi'\cdot\phi & 2\psi''\cdot\phi \\ \psi''\cdot\phi & -1 & \psi'\cdot\phi \\ \psi'\cdot\phi & 0 & -1 \end{vmatrix} = \mathbf{D}_y^3\phi$$

ammettendo che sviluppato il determinante si considerino gli esponenti di $\phi$ quali indici di differenziazione. In generale si otterà.

$$(-1)^n \mathbf{D}_y^{n+1}\phi =$$

$$\begin{vmatrix} \psi^{(n+1)}\phi & \psi'\phi & n\psi''\phi & \frac{n(n-1)}{1\,.\,2}\psi'''\phi & \ldots & \frac{n(n-1)}{1\,.\,2}\psi^{(n-1)}\phi & n\psi^{(n)}\phi \\ \psi^{(n)}\phi & -1 & \psi'\phi & (n-1)\psi''\phi & \ldots & \frac{(n-1)(n-2)}{1\,.\,2}\psi^{(n-2)}\phi & (n-1)\psi^{(n-1)}\phi \\ \psi^{(n-1)}\phi & 0 & -1 & \psi'\phi & \ldots & \frac{(n-2)(n-3)}{1\,.\,2}\psi^{(n-3)}\phi & (n-2)\psi^{(n-2)}\phi \\ \ldots & \ldots & \ldots & \ldots & \ldots & \ldots & \ldots \\ \ldots & \ldots & \ldots & \ldots & \ldots & \ldots & \ldots \\ \ldots & \ldots & \ldots & \ldots & \ldots & \ldots & \ldots \\ \psi''\phi & 0 & 0 & 0 & \cdots & -1 & \psi'\phi \\ \psi'\phi & 0 & 0 & 0 & \cdots & 0 & -1 \end{vmatrix}$$

Parigi 26 dicembre 1855.

---